\documentclass[12pt]{article}
\usepackage{amsmath} 
\usepackage{amssymb} 
\usepackage{latexsym}
\usepackage{theorem} 
\def\mylabel#1{\label{#1}}

\newtheorem{theorem}{Theorem}
\newtheorem{lemma}[theorem]{Lemma}
\newtheorem{corollary}[theorem]{Corollary}
\newtheorem{proposition}[theorem]{Proposition}

{\theorembodyfont{\rmfamily}

}
{\theoremstyle{break}	
{\theorembodyfont{\rmfamily} 
				\newtheorem{example}[theorem]{Example}
				
}}
{\end{enumerate}}

\newenvironment{proof}[1]{\smallskip \noindent {\bf #1}}{\qed\smallskip}

\def\qed{\ifhmode\unskip\nobreak\fi\ifmmode\ifinner\else\hskip5pt\fi\fi
 \hfill\hbox{\hskip5pt\vrule width4pt height6pt depth1.5pt\hskip1pt}}

\newcommand{\ov}[1]{\ensuremath{\overline{#1}}} 
\newcommand{\RR}{\ensuremath{{\bf R}}}     
\newcommand{\R}[1]{\ensuremath{{\bf R}^{#1}}} 
\newcommand{\CC}{\ensuremath{{\bf C}}}     

\newcommand{\Fix}[1]{\ensuremath{\mathsf {Fix}\,(#1)}}
\newcommand{\co}{\colon\thinspace} 

\newcommand{\ca}[1]{\ensuremath{{\cal #1}}}

\begin{document}

\title{\bf Fixed points of analytic actions of supersoluble Lie groups
on compact surfaces}

\author
{
{\bf Morris W. Hirsch\thanks{Partially supported by NSF grant DMS-9802182}~
\&  Alan Weinstein\thanks{Partially supported by NSF Grant DMS-9971505}
}
\\ Department of
Mathematics\\University of California at Berkeley
}
\maketitle

\begin{abstract} We show that every real analytic action of a connected
supersoluble Lie group on a compact surface with nonzero Euler
characteristic has a fixed point.  This implies that E. Lima's fixed
point free $C^{\infty}$ action on $S^2$ of the affine group of the
line cannot be approximated by analytic actions.  An example is given
of an analytic, fixed point free action on $S^2$ of a solvable group
that is not supersoluble.
\end{abstract}

\section*{Introduction}   \mylabel{sec:intro}

Let $M$ denote a compact connected surface, with possibly empty
boundary $\partial M$, endowed with a (real) analytic
structure. $T_pM$ is the tangent space to $M$ at $p\in M$.  The Euler
characteristic of $M$ is denoted by $\chi (M)$.

Let $G$ be a Lie group with Lie algebra $\ca L (G)=\ca G$; all groups
are assumed connected unless the contrary is indicated.  
An {\em action} of $G$ on $M$ is a homomorphism $\alpha$ from $G$ to
the group $\mathsf{H} (M)$ of homeomorphisms of $M$ such that the
evaluation map
\[\mathsf{ev}^\alpha=\mathsf{ev}\co G\times M\to M,\;(g,x)\to\alpha (g)(x)\]
is continuous.  We usually suppress notation for $\alpha$, denoting
$\alpha (g)(x)$ by $g(x)$.  The action is called $C^r,\:r\in
\{1,2\dots;\omega\}$ if $\mathsf{ev}$ is a $C^r$ map, where $C^\omega$
means analytic.

The set $\ca A (G, M)$ of actions of $G$ on $M$ is embedded in the
space of continuous maps $G\times M\to M$ by the correspondence
$\alpha\mapsto\mathsf{ev}^\alpha$.  We endow $\ca A (G, M)$ with the
topology of uniform convergence on compact sets.

A point $p\in M$ is a {\em fixed point} for an action $\alpha$ of $G$
if $\alpha(g)(p)=p$ for all $g\in G$.  The set of fixed points is
denoted by $\Fix G$ or $\Fix {\alpha(G)}$.

In this paper we consider the problem of finding conditions on
solvable group actions that guarantee existence of a fixed point.

When $\chi(M)\ne 0$, every flow (action of the real line $\RR$) on $M$
has a fixed point; this was known to Poincar\'e for flows generated by
vector fields, and for continuous actions it is a well known
consequence of Lefschetz's fixed point theorem.  E. Lima \cite{Li64}
showed that every abelian group action on $M$ has a fixed point, and
J. Plante \cite{Pl86} extended this to nilpotent groups.

These results do not extend to solvable groups: Lima \cite{Li64}
constructed a fixed point free action on the 2-sphere of the solvable
group $A$ of homeomorphisms of $\RR$ having the form $x\mapsto ax+b,
\;a>0, b\in\RR$; and Plante \cite{Pl86} constructs fixed point free
action of $A$ on all compact surfaces.  These actions are not known to
be analytic; but Example \ref{th:ex1} below describes a fixed point
free, analytic action of a 3-dimensional solvable group on $S^2$.

Recall that $G$ is {\em supersoluble} if 
 every element of $\ca
G$ belongs to a codimension one subalgebra  (see Barnes \cite{Ba67}).
Our main result is the following theorem:
\begin{theorem}		\mylabel{th:main}
Let $G$ be a connected supersoluble Lie group and $M$ a compact
surface $M$ such that $\chi (M)\ne 0$.  Then every analytic action of
$G$ on $M$ has a fixed point.
\end{theorem}

Since the group $A$ described above is supersoluble, Lima's
$C^\infty$ action cannot be improved to a fixed point free analytic
action.  The following result shows it cannot be approximated by
analytic actions:

\begin{corollary}		\mylabel{th:maincor}
Let $G$ and $M$ be as in Theorem \ref{th:main}.  If $\alpha \in \ca A
(G, M)$ has no fixed point, then $\alpha$ has a neighborhood in $\ca A
(G, M)$ containing no analytic action.
\end{corollary}
\begin{proof}
By Theorem \ref{th:main} and compactness of $M$, it suffices to prove
the following: For all convergent sequences $\beta_n\to \beta$ in $\ca A
(G,M)$ and $p_n \to p$ in $M$, with $p_n\in \Fix{\beta_n(G)}$, we have
$p\in\Fix{\beta (G)}$.  Being a connected locally compact group, $G$
is generated by a compact neighborbood $K$ of the identity.  Then
$\beta_n (g)\to \beta (g)$ uniformly for $g\in K$, so $\beta (g)
(p)=p$ for all $g\in K$.  Since $K$ generates $G$, this implies that
$p\in\Fix{\beta (G)}$.
\end{proof}

In Theorem \ref{th:main}, the hypothesis that $G$ is connected is
essential: the abelian group of rotations of $S^2$ generated by
reflections in the three coordinate axes is a well known
counterexample.  And every Lie group with a nontrivial homomorphism to
the group of integers acts analytically without fixed point on every
compact surface admitting a fixed point free homeomorphism, thus on
every surface except the disk and the projective plane.

The following example shows that supersolubility is essential:

\begin{example}		\mylabel{th:ex1}
Let $Q$ be the 3-dimensional Lie group obtained as the semidirect
product of the real numbers $\RR$ acting on the complex numbers $\CC$
by $t\cdot z=e^{it}z$; this group is solvable but not
supersoluble.  Identify $Q$ with the space $\RR\times \CC\approx\R 3$
and note that left multiplication defines a linear action of $Q$ on
$\R 3$.  The induced action on the 2-sphere $S$ of oriented
lines in $\R 3$ through the origin has no fixed point, and $\chi
(S)=2$.  Geometrically, one can see this as the 
universal cover of the proper euclidean motions of the plane, acting
on two copies of the plane joined along a circle at infinity.  
\end{example}

\bigskip
We thank F.-J. Turiel for pointing out a small error in
an earlier version of our manuscript.  He has also obtained some
interesting results complementary to ours in \cite{tu-analytic}.

\section*{Proof of Theorem \ref{th:main}}
We assume given an action $\alpha\co G\to \mathsf{H} (M)$.
The {\em orbit} of $p\in M$  is $G(p)= \{ g(x) \co g \in G \}$.
The {\em isotropy group} of $p\in M$ is the closed subgroup
$I_p=\{g\in G\co \alpha(g)(p)=p\}$.
The {\em evaluation map} $\mathsf{ev}_p\co G\to M$ at $p\in M$ is
defined by $g\mapsto g(p)$. 

Suppose that the action is $C^r, \:r\ge 1$.  Then $\mathsf{ev}_p$ induces a
bijective $C^r$ immersion $i_p\co G/I(p) \to G(p)$.  The tangent space
$E(p)\subset T_p M$ to this immersed manifold at $p$ is the image of
$T_e G$ under the differential of $\mathsf{ev}_p$ at the identity
$e\in G$.

For $j=0,1,2$, let $V_j=V_j(G)\subset M$ denote the union of the
$j$-dimensional orbits.  Then $M=V_2\cup V_1\cup V_0$. Each $V_j$ is
invariant, $V_2$ is open, $V_1\cup V_0$ is compact, and $V_0=\Fix G$.

\begin{lemma}[Plante]		\mylabel{th:plante}
Assume that $G$ is solvable and that $G(p)$ is a compact 1-dim\-en\-sion\-al orbit.
Then there is a closed normal subgroup $H\subset G$ of codimension 1
such that every point of $G(p)$ has isotropy group $H$.
\end{lemma}
\begin{proof}
Choose a homeomorphism $f\co G(p)\approx S^1$ (the circle).  Let
$\beta\co G\to \mathsf{H}(S^1)$ be the action defined by $\beta(g) =
f\circ \alpha (g)\circ f^{-1}$.  Because $G$ is solvable, by a result
of Plante (\cite{Pl86}, Theorem 1.2) there exists a homeomorphism $h$
of $S^1$ conjugating $\beta(G)$ to the rotation group
$\mathrm{SO}(2)$.  Since $\beta(G)$ is abelian and acts transitively
on $S^1$, all points of $S^1$ have the same isotropy group for
$\beta$; this isotropy group is the required $H$.
\end{proof}

Analyticity is used to establish the following useful property:

\begin{lemma}		\mylabel{th:lem2}
Assume that $G$ acts analytically and that $\Fix G=\varnothing$.  Then either
$V_1=M$ and $\chi (M)=0$, or else $V_1$ is the (possibly empty) union of
a finite family of orbits, each of which is a smooth Jordan curve
contained in $\partial M$ or in $M\setminus \partial M$.
\end{lemma}
\begin{proof}
Since there are no orbits of dimension $0$, $V_1$ is a compact set
comprising the  points $p$ such that $\dim E_p \le 1$. It is easy
to see that $V_1$ is a local analytic variety.

If $V_1=M$ then the map $p\mapsto E_p$ is a continuous field of tangent
lines to $M$, tangent to $\partial M$ at boundary points.  The
existence of such a field implies that $\chi (M)= 0$.

Assume that $V_1\ne M$.   Note that $\dim_p V_1\ge 1$ at each $p\in V_1$.
Since $M$ is connected and $V_1$ is a variety, $V_1$ must have dimension
$1$ at each point.  The set of points where $V_1$ is not smooth is a
compact, invariant 0-dimensional subvariety, i.e., a finite set of
fixed points, hence empty.  Since $V_1$ consists of 1-dimensional
orbits, $V_1$ must be a compact, smooth
invariant 1-manifold without boundary, i.e.
each component of $V_1$ is a Jordan curve.  
Since
$\partial M$ is the union of invariant Jordan curves, any component of
$V_1$ that meets $\partial M$ is a component of $\partial M$.
\end{proof}

In view of Lemma \ref{th:lem2}, it suffices to prove the following
more general result:

\begin{proposition}		\mylabel{th:main2}
Let $G$ be a connected supersoluble Lie group acting continuously on
the compact connected surface $M$.  Assume that
\begin{description}
\item [(a)] there are no fixed points
\item [(b)] for each closed subgroup $H,\,$ $V_1(H)$ is the union
(perhaps empty) of finitely many disjoint Jordan curves. 
\end{description}
 Then $\chi
(M)=0$.
\end{proposition}
By passing to a universal covering group we assume that $G$ is simply
connected.  This implies that every closed subgroup is simply
connected (see Hochschild \cite{Ho65}, Theorem XII.2.2.)
 
We proceed by induction on $\dim G$, the case $G=\RR$ having been
covered in the introduction.  Henceforth assume inductively that $\dim
G = n\ge 2$ and that the proposition holds for all supersoluble groups
of lower dimension.  With this hypothesis in force, we first rule out
the case that $M$ is a disk:

\begin{proposition}		\mylabel{th:propD}
If $M$ is as in Proposition \ref{th:main2}, then $\chi (M)\ne 1$
\end{proposition}
\begin{proof}
Suppose not; then $M$ is a closed 2-cell.  Since there are no fixed
points, $\partial M$ is an orbit, hence a component of $V_1$.  Every
component of $V_1$ bounds a unique 2-cell in $M$, and there are only
finitely many such 2-cells.  Let $D$ be one that contains no other.
Then $D$ is invariant under $G$, and the action of $G$ on $D$
is fixed point free.    Therefore we may assume
that $M=D$, so that $V_1=\partial M$.   

By Lemma \ref{th:plante} there exists a closed normal subgroup $H$ of
codimension one with $\partial M\subset \Fix H$.  Let $R\subset G$ be
a 1-parameter subgroup transverse to $H$ at the identity; then
$RH=G$. 

Because $G$ is supersoluble, there is a codimension one subalgebra
$\ca K\subset \ca G$ containing the Lie algebra $\ca R$ of $R$.  Because $G$ is simply
connected and solvable $\ca K$ is the Lie algebra of a closed subgroup
$K\subset G$ of dimension $n-1$, and $KH=G$.
By the induction hypothesis there exists $p\in\Fix K$. Then $\dim G(p)
\le \dim G -\dim K =1$.  Therefore $p\in V_1=\partial D$.  We now have
$p\in\Fix K\cap\Fix H=\Fix G$, a contradiction.
\end{proof}

We return now to the case of general $M$.

Denote the connected components of $M\setminus V_1$ by
$U_i,\dots,U_r,\; r\ge 1$.  Each $U_i$ is an open orbit, whose set
theoretic boundary $\mathsf {bd}\,U_i$ is a (possibly empty) union of
components of $V_1$.  The closure $\ov U_i$ is a compact surface
invariant under $G$, whose boundary as a surface is $\partial
U_i=\mathsf {bd}\,U_i$.

We show that $U_i$ is an open annulus.  Let $H\subset G$ be the
isotropy subgroup of $p\in U_i$.  Evaluation at $p$ is a surjective
fibre bundle projection $G\to U_i$ with standard fibre $H$.  Therefore
there is an exact sequence of homotopy groups
\[\cdots\to\pi_j (G)\to\pi_j (U_i)\to \pi_{j-1}(H)\to \pi_{j-1} (G)\to\dots\to
\pi_0 (G)=\{0\}
\]   
ending with the trivial group $\pi_0 (G)$ of components of $G$.
The component group $\pi_0(H)$ is solvable (see Raghunathan \cite{Ra72},
Proposition III.3.10), so taking $j=1$ shows that $\pi_1 (U_i)$ is
solvable. 
Therefore $U_i$ is a
sphere, torus, open 2-cell, or open annulus.  If $U_i$ is a torus then
$U_i=M$, contradicting $\chi (M)\ne 0$.  The sphere is ruled out by
the exact sequence $\pi_2(G)\to\pi_2(U_i)\to\pi_1(H)$, because $\pi_2
(G)={0}$ for every Lie group and $\pi_1(H)={0}$. 
Proposition \ref{th:propD} rules out the 2-cell.

It follows that $\ov U_i$ is a closed annulus, so $\chi (\ov U_i)=0$.
By the additivity property $\chi(A\cup B) = \chi(A) + \chi(B) -
\chi(A\cap B)$ of the Euler characteristic, any space $M$ built by
gluing annuli along their boundary circles must have $\chi(M)=0$.

\qed


\end{document}